\documentclass[11pt]{article}
\usepackage{comment}
\usepackage{tabularx,lscape,longtable}
\usepackage{comment}
\usepackage{amsfonts}
\usepackage{graphicx}
\usepackage{epsfig}
\usepackage{epstopdf}
\usepackage{multirow}
\usepackage{latexsym}
\usepackage{amssymb}
\usepackage{amsmath}
\usepackage{array}
\usepackage{tikz}
\usepackage{tfrupee}
\usepackage{xcolor}
\usepackage{multicol}
\definecolor{cadmiumred}{rgb}{10,0,2}
\definecolor{red}{rgb}{1,0,0}
\definecolor{blue}{rgb}{0,0,1}
\definecolor{green}{rgb}{0,1,0}
\begin{document}
\newtheorem{t1}{Theorem}[section]
\newtheorem{d1}{Definition}[section]
\newtheorem{c1}{Corollary}[section]
\newtheorem{l1}{Lemma}[section]
\newtheorem{r1}{Remark}[section]
\newtheorem{ex}{Example}[section]
\newtheorem{re1}{Result}[section]
\newtheorem{co}{Counterexample}[section]
\newtheorem{p1}{Proposition}[section]
\title{On weighted failure rate, its means and  associated quantile version}
 \author{Subarna Bhattacharjee$^{1}$$^{}$\thanks{Corresponding author:~~
E-mail: subarna.bhatt@gmail.com}, S.M. Sunoj$^{2}$, Sabana Anwar$^{3}$ \\
{\it $^{1,3}$ Department of Mathematics, Ravenshaw University, Cuttack, Odisha, India} \\
{\it $^{2}$ Department of Statistics, Cochin University of Science and Technology, Cochin, Kerala, India}
\\}
\date{\today} 
\maketitle
\begin{abstract}
 In this paper, we define weighted failure rate and their different means from the stand point of an application. We begin by emphasizing that the formation of $n$ independent component series system having  weighted failure rates with sum of weight functions being unity is same as a mixture of $n$ distributions. We derive some parametric and non-parametric characterization results. We discuss on the form invariance property of baseline failure rate for a specific choice of weight function. Some bounds on means of aging functions are obtained.  Here, we establish that weighted IFRA class is not closed under formation of coherent systems unlike the IFRA class. An interesting application of the present work is credited to the fact that the quantile version of means of failure rate is obtained as a special case of weighted means of failure rate.
\end{abstract}
\section{Introduction}
The notion of ageing plays an important role in reliability theory and in the study of lifetime data analysis.  Ageing of a mechanical or biological component based on a lifetime distributions is generally studied using the residual lifetime of the unit that is affected its age.  Abundant literature is available on various ageing concepts and their patterns of ageing, comparison of life distributions and to explain their data generating mechanism.   Reliability ageing classes based on the monotonicity of the failure rate, such as increasing (decreasing) failure rate (IFR (DFR)) and its average, increasing (decreasing) failure rate average (IFRA (DFRA)) have been found great interest among researchers as it easily give an indication on the manner in which ageing can be described, life distributions can be classified and distinguished, and appropriate models can be chosen when observations are available (cf. Barlow and Proschan (1975)). \\

Let $X$ be a non-negative random variable representing the lifetime of an event or living mechanism with absolutely continuous distribution function $F(\cdot)$ and probability density function (pdf), $f(\cdot)$.  Then $F$ is said to be IFR (DFR), if the conditional survival function $\Bar{F}(x|t) = \frac{\Bar{F}(x + t)}{\Bar{F}(x)}$ is decreasing (increasing) in $0 \leq t < \infty$, where $\Bar{F} = 1 - F$ is the survival (reliability) function; or equivalently the failure rate $h(t) = \frac{f(t)}{\Bar{F}(t)}$ is increasing (decreasing) in $t \geq 0$, provided $f(t)$ exists.  Further, $F$ is said to IFRA (DFRA), if $-\left(\frac{1}{t}\right)\log \Bar{F}(t)$ is increasing (decreasing) in $t \geq 0$.  However, in many real situations, $h(x)$ is not always monotonic.  In such cases, the monotonicity of IFRA class condition in terms of the failure rate, $\frac{1}{t} \int_{0}^{x}h(t)dt$, known as the arithmetic mean failure rate (AFR) is a useful measure (Roy and Mukherjee (1992)) in identifying the monotonicity of classes of life distributions.  Along with arithmetic mean failure rate, Roy and Mukherjee (1992) have also studied classes of distributions through the monotonic behaviour of geometric failure rate (GFR) and harmonic failure rate (HFR), and the characterizations and ageing classes based on it.
They pointed out until then no work has been done on GFR and HFR.
The following definition is cited from Roy and Mukherjee (1992).
\begin{d1}\label{d1.1}
Let $X$ be a non-negative random variable with absolutely continuous CDF $F(\cdot)$, PDF $f(\cdot)$ and failure rate $h(\cdot)$. Then the arithmetic mean failure rate (AFR), geometric mean failure rate (GFR) and harmonic failure rate (HFR), denoted by $A(\cdot), G(\cdot)$ and $H(\cdot)$ respectively are defined as
$A(x)=\frac{1}{x}\int_{0}^{x}h(u)du,x>0;$ $G(x)=\exp\Big(\frac{1}{x}\int_{0}^{x}\ln h(u)du\Big)$;
 $H(x)=\Big(\frac{1}{x}\int_{0}^{x}\frac{1}{h(u)}du\Big)^{-1}.$$\hfill\Box$
\end{d1}
%A note on applications of different formulae on these means can %be referred from Bullen (2003).  
Recently, Bhattacharjee et al. (2022) further studied the usefulness of the three measures based on the notion of ageing intensity (AI) function proposed by Jiang et al. (2003).\\

When sample observations are not equally likely, we use the weighted measures to capture the significance of their relative importance.  Choosing appropriate weights, we can then compute various measures in a better way by giving appropriate weights based on the sample mechanism.  Such biased sampling schemes are usually employed in observational studies either due to 
its convenience or its cost-effectiveness.  Based on this, Rao (1965) identified the concept of weighted distributions in connection with the modeling statistical data, in situations where the usual practice of employing standard distributions for the purpose was not found appropriate.  These distributions occur frequently in the studies related to reliability, analysis of family data, meta analysis and analysis of intervention data, biomedicine, ecology, etc, for more details, see Patil and Rao (1978), Gupta and Kirmani (1990),  and the references therein.  If $X$ is a non-negative random variable with a probability density function (pdf) $f(x)$, then the pdf of the weighted random variable $X^w$ is given by, $f_w (x) = \frac{w(x) f(x)}{Ew(X)}, x > 0$, where $w(\cdot)$ is a non-negative weight function (cf. Rao (1965)).  There are many weight functions used by different authors, however, the weight functions $w(x) = x$ and $w(x) = x^{\alpha}, \alpha > 0$ are found to be more popular due to its adaptibility in terms of identifying the observed distribution in various applied problems wherein the probability of selecting the sample units are proportional to the length or size of the population units, the respective random variables are known as the length-biased and size-biased random variables.  Motivated by these, in the present study, we propose weighted mean failure rates based on the measures of AFR, GFR and HFR.\\

\begin{comment}
A distribution can be specified either using the distribution function or its quantile function. The quantile funcion of a random variable $X$ is defined as
\begin{equation}\label{1.1}
    Q(u) = F^{-1}(u) = \inf \{x: F(x) \geq u \}, \; 0 \leq u \leq 1.
\end{equation}
From \eqref{1.1}, we have $F\left(Q(u) \right) = u$, differentiating it with respect to $u$, we get $f\left( Q(u)\right) q(u) = 1$, where $f\left( Q(u)\right)$ and $q(u) = \frac{dQ(u)}{du}$ are respectively known as the density quantile function and quantile density function of the random variable $X$.
\end{comment}
The paper is organized as follows. In Section 2, we introduce the definitions of weighted means of failure rate along with citing an application. As an application of weighted failure rates as proposed in this paper, we note that formation of an $n$ component series system having complementary weight functions is actually a mixture of $n$ distributions and vice-versa. We give some parametric and non-parametric characterization results. Section 2 gives a note on form invariance property of the baseline failure rate and its transformation from one aging class to another depending upon the choice of the weight function. Section 2 also  highlights on some bounds of means of aging functions. We focus on some new non-parametric aging classes based on means of failure rate and discuss their inclusive properties. Some illustrative examples are given for ready reference. Some equivalent conditions of aging classes based on geometric and harmonic means are obtained. We prove our claim that weighted IFRA class is not closed under formation of coherent systems unlike IFRA class by giving an easy counterexample. In Section 4, we derive the quantile version of means of failure rate as special case from weighted means of failure rate. We study the proportional quantile hazards model and compare it with conventional proportional hazards model. Concluding remarks are listed in Section 5.

\section{Weighted means of failure rate}

In this paper, we introduce a generalized versions of AFR, GFR and HFR involving a suitable choice of a non-negative weight function as defined below.
\begin{d1}
\label{newd}
Let $X$ be a non-negative random variable with absolutely continuous distribution function $F(\cdot)$, probability density function $f(\cdot)$ and failure rate $h(\cdot).$ The weighted arithmetic mean failure rate (w-AFR), weighted geometric mean failure rate (w-GFR) and weighted harmonic failure rate (w-HFR) denoted by $A^{w}(\cdot), G^{w}(\cdot)$ and $H^{w}(\cdot)$ respectively, with a suitable non-negative weight function $w(\cdot)$, are defined as
\begin{enumerate}
\item[(i)] $A^{w}(x)=\frac{\int_{0}^{x}w(u)h(u)du}{\int_{0}^{x}w(u)du}, x>0;$
\item[(ii)] $G^{w}(x)=\exp\Big(\frac{\int_{0}^{x}w(u)\ln h(u)du}{\int_{0}^{x}w(u)du}\Big), x>0;$
\item[(iii)]$H^{w}(x)=\Big(\int_{0}^{x}w(u)du\Big)\Big(\int_{0}^{x}\frac{w(u)}{h(u)}du\Big)^{-1}, x>0.$$\hfill\Box$
\end{enumerate}
\end{d1}
Clearly, if $w(x)=1$ for all $x\geq 0$ then above definition reduces to that of AFR, GFR and HFR given in Definition \ref{d1.1} due to Roy and Mukherjee (1992).\\ 

In the pretext, of above, we shall define the other reliability functions of the weighted random variable as given in the following definition.
\begin{d1}
\label{newd}
The weighted survival function of $X,$ or survival function of weighted random variable $X^{w},$ denoted by $\overline{F}^{w}(\cdot)$ is defined as $\overline{F}^{w}(x)=\exp\big(-\int_{0}^{x}w(u)h(u)du\big), x>0.$ The density and failure rate function of $X^{w}$ are $f^{w}(x)=w(x)h(x)\exp\big(-\int_{0}^{x}w(u)h(u)du\big)$ and $h^{w}(x)=w(x)h(x)$ for all $x>0,$ respectively.
\end{d1}
The fact that $h^{w}(x)=w(x)h(x)$ reminds us of proportional hazard models (PHR) where $h(x)$ is the baseline failure rate and $w(x)$ is the proportionality function giving rise to a new hazard rate $h^{w}(\cdot).$ \\

Referring to the related literature, one can notice that corresponding to  the baseline survival function $\overline{G}$ having failure rate $r_{G},$ Marshall and Olkin (1997) proposed a cumulative distribution function $F$ such that its hazard rate $h_{F}(\cdot)$ is given by $h_{F}(x,\alpha)=\frac{1}{1-\overline{\alpha}\overline{G}(x)}h_{G}(x)$ where $x,\alpha\in R^{+}$ and $\overline{\alpha}=1-\alpha,$ ( the parameter $\alpha$ termed as tilt parameter by Marshall and Olkin (2007)) and this is a special case of Definition \ref{newd} if one assumes $w(x)=\frac{1}{1-\overline{\alpha}\overline{G}(x)}.$ Furthermore, Balakrishnan et al. 
(2018) defined modified proportional hazard rates (MPHR) of $n$ independent components having lifetimes $X_{1}, X_{2}, \ldots, X_{n}$ with respective survival functions $\overline {F}_{i}$ if  $$F_{i}(x,\lambda_{i})=\frac{1-\big(\overline{F}(x)\big)^{\lambda_{i}}}{1-\overline{\alpha}\big(\overline{F}(x)\big)^{\lambda_{i}}},\alpha>0, \overline{\alpha}=1-\alpha, \lambda_{i}>0$$ \mbox{~for~} $i=1,2,\ldots,n$ where $\overline{F}$ is the corresponding baseline survival function. They considered it (MPHR model) to be the generalization of PHR model because if $\alpha=1$ then PHR is a special case of MPHR. 
However, one shall observe that this is based on the notion that $X_{i}$'s with survival functions $\overline{F}_{i}(x)$ follow PHR model if there exits positive constants $\lambda_{i}$'s such that $\overline{F}_{i}(x)=\overline{F}^{\lambda_{i}}(x).$  It is worthwhile to note that the definition of MPHR proposed by Balakrishnan et al. 
(2018) is not based on the fact that ${h}_{i}(x)=h(x)w(x)$ where $h(x)$ is considered to be the baseline failure rate. The present work, in other words, is an attempt to define PHR model in a more general sense. \\

We now give a necessary and sufficient condition that  weight function $w(x)$, and hazard rate $h(x)$ must satisfy so that $\bar{F}^{w}(x)$ represents a (weighted) survival function. One can refer to Marshall and Olkin (2007) to look into the postulates for hazard rate (non-weighted).
\begin{enumerate}
\item[(i)] $w(x)\geq 0,$ $h(x)\geq 0.$
\item[(ii)] For $x>0,$ $\int_{0}^{x} w(u)h(u) du <\infty$
\item[(iii)]$\int_{0}^{\infty} w(u)h(u) du =\infty$
\item[(iv)] If $\int_{0}^{x} w(u)h(u) du =\infty$ for some $x$ then $h(y)=\infty$ for every $y>x.$
\end{enumerate}

Now, we look into some uses of weighted failure rate arising in practical field. Let us consider a series system formed by $n$ components having failure rates $h_{i}(x)$ with respective weights $w_i(x),$  for $i=1,2,\ldots,n$ and $x>0$ such that $\sum_{i=1}^{n}w_{i}(x)=1.$ The failure rate $h(x)$ of the resultant $n$ component series system is $h(x)=\sum_{i=1}^{n}h_{i}(x)w_{i}(x).$ This form of $h(x)$ is similar to the failure rate of the mixture of $n$ distributions, with cumulative distributions, $F_{i}(\cdot)$  having failure rates, $h_{i}(\cdot)$ for $i=1,2,\ldots,n.$ The failure rate of mixture of $n$ distributions is given by $F(x)=\sum_{i=1}^{n}\pi_{i} F_{i}(x),$ with $\sum_{i=1}^{n}\pi_{i}=1$ is 
\begin{eqnarray}
h(x)&=&\frac{\sum_{i=1}^{n}\pi_{i} f_{i}(x)}{1-\sum_{i=1}^{n}\pi_{i} {F}_{i}(x)}\nonumber\\
&=&\frac{\sum_{i=1}^{n}\pi_{i}h_{i}(x)\overline{F}_{i}(x)}{\sum_{i=1}^{n}\pi_{i}\overline{F}_{i}(x)}
=\sum_{i=1}^{n}p_{i}(x)h_{i}(x)\nonumber
\end{eqnarray}
where $p_{i}(x)=\frac{\pi_{i} \overline{F}_{i}(x)}{\sum_{i=1}^{n}\pi_{i} \overline{F}_{i}(x)}$ which in turn satisfies $\sum_{i=1}^{n}p_{i}(x)=1.$\\ 

Unlike the log-exponential family (cf. Patil G. P.  and  Ord, J. K. (1976)) which possesses the  form-invariance property among the  weighted distributions defined by Rao (1965) under size based sampling of order $c >0$, (i.e., $w(x)=x^{c}$),  in the present work Weibull distribution bear the said property. If a given Weibull distribution with failure rate $h(x)=\alpha\beta x^{\beta-1}$ belongs to positive aging classes, namely IFR, then the resultant size biased  distribution also fall in the same aging class, but this may not be true for DFR negative aging class, i.e., if $\beta+c>1$ then the baseline DFR Weibull distribution is shifted to IFR positive aging class under size based sampling.\\

Additive Weibull distribution with failure rate $h(x)=\alpha \theta x^{\theta-1}+\beta\gamma x^{\gamma-1}$ with $\alpha, \beta,\theta,\gamma >0$ has form-invariance property under size biased sampling. However, the weight function  $w(x)=x^{c}$ drags the additive Weibull distribution from DFR class to IFR class if $c+\theta >$ and $c+\gamma >1.$ \\

If $X$ follows four-parameter Weibull distribution (cf. Kies (1958)) with survival function  $$\overline{F}(t)=\exp\big(-\lambda \big(\frac{t-a}{b-t}\big)^{\beta}\big), 0\leq a<t<b, \lambda, \beta>0$$ then $w(t)=\big(\frac{t-a}{b-t}\big)$ is a form-invariance weight function for the distribution. We know that $X$ has bathtub failure rate if $0<\beta<1,$ and IFR if $\beta>1$ but the weighted random variable $X^{w}$ is always IFR independent of $\beta$ under the aforementioned weighted transformation. 
\begin{comment}
  \begin{ex}
Let $X$ follows two-parameter Weibull distribution with scale parameter $\alpha$ and $\beta.$ 
If we consider $w(t)=e^{nt}$ and $h(t)=\alpha \beta t^{\beta -1}$  then $\bar{F}(t)=\exp\Big(-t^\beta (-nt)^{-\beta} \alpha \beta~\Big(\Gamma[\beta] - \Gamma[\beta, -n t]\Big)\Big)$
\end{ex}
\end{comment}
\subsection{Characterization results}
Now, we prove that the equality of any two of w-AFR, w-GFR and w-HFR characterize exponential distribution. The following proposition gives some light on this. We continue with the same notations as discussed in previous sections of the paper.
\begin{p1}
\label{pp}
A non-negative random variable $X,$ follows exponential distribution if and only if for $x>0,$ any one of the following hold\begin{enumerate}
\item[(i)] $A^{w}(x)=G^{w}(x)$
\item[(ii)] $G^{w}(x)=H^{w}(x)$
\item[(iii)]$A^{w}(x)=H^{w}(x)$.
\end{enumerate}
\end{p1}
{\bf Proof.} If $X$ follows exponential distribution then it is easy to prove that (i), (ii) and (iii) hold. Conversely, if (i) holds then $$\frac{\int_{0}^{x}w(u)h(u)du}{\int_{0}^{x}w(u)du}=\exp\Big(\frac{\int_{0}^{x}w(u)\ln h(u)du}{\int_{0}^{x}w(u)du}\Big),$$ gives
\begin{equation}
\label{wt1} {\scriptsize \Big(\int_{0}^{x}w(u)du\Big)\Big\{\ln\Big(\int_{0}^{x}h(u)w(u)du\Big)\Big\}=\Big(\int_{0}^{x}w(u)du\Big)\Big(\ln \int_{0}^{x}w(u)du\Big)+\int_{0}^{x}w(u)\ln h(u)du}.
\end{equation}
Differentiating (\ref{wt1}) with respect to $x,$ we get $\ln (ez_1(x))=z_1(x)$ where $$z_{1}(x)=\frac{h(x)\int_{0}^{x}w(u)du}{\int_{0}^{x}w(u)\ln w(u) h(u)du}.$$ Thus, $z_1(x)=1,$ for all $x \geq 0,$ gives $\Big(\frac{d}{dx}h(x)\Big(\int_{0}^{x}w(u)du\Big)=0,$ and since $\int_{0}^{x}w(u)du\neq 0$, we conclude that $h(x)$ is constant for all $x\geq 0.$ This proves that $X$ has exponential distribution. Similarly, if $(ii)$ holds then $$\exp\Big(\frac{\int_{0}^{x}w(u)\ln h(u)du}{\int_{0}^{x}w(u)du}\Big)=\Big(\int_{0}^{x}w(u)du\Big)\Big(\int_{0}^{x}\frac{w(u)}{h(u)}du\Big)^{-1},$$ 
or equivalently
\begin{equation}
\label{wt2}
\int_{0}^{x}w(u)h(u)du+\Big(\int_{0}^{x}w(u)du\Big)\ln \Big(\int_{0}^{x}\frac{w(u)}{r(u)}du\Big)=\Big(\int_{0}^{x}w(u)du\Big)\ln\Big(\int_{0}^{x} w(u)du\Big).
\end{equation}
After differentiating (\ref{wt2}), we get $\ln (ez_{2}(x))=z_{2}(x)$ where $$z_{2}(x)=\frac{\int_{0}^{x}w(u)du}{h(x)\Big(\int_{0}^{x}\frac{w(u)}{h(u)}du\Big)}.$$
Hence, $z_{2}(x)=1$ which in turn gives $\frac{d}{dx}r(x)\Big(\int_{0}^{x}\frac{w(u)}{h(u)}du\Big)=0.$ Since $\Big(\frac{w(u)}{\int_{0}^{x}h(u)}du\Big)\neq 0,$ it follows that $h(x)=\mbox{constant}.$ This proves that if $(ii)$ holds then $X$ has exponential distribution.  \\

Note that if $(iii)$ holds then it is equivalent to the fact that 
\begin{equation}
\label{wt3}
\Big(\int_{0}^{x}w(u)h(u)du\Big)\Big(\int_{0}^{x}\frac{w(v)}{h(v)}dv\Big)=\Big(\int_{0}^{x}w(u)du\Big)^{2}.
\end{equation}
Taking logarithm on both sides and then differentiating both sides with respect to $x,$ we get
\begin{equation}
\label{cha}
\frac{h(x)\int_{0}^{x}w(u)du}{\int_{0}^{x}w(u)h(u)du}+\Big(\frac{1}{h(x)}\Big)\Big(\frac{\int_{0}^{x}w(u)du}{\int_{0}^{x}\frac{w(v)}{h(v)}dv}\Big)=2.
\end{equation}
Since w-HFR=w-AFR,  replacing w-HFR by w-AFR in the second term of  (\ref{cha}), we get
\begin{displaymath}
\label{cha1}
\frac{h(x)\int_{0}^{x}w(u)du}{\int_{0}^{x}w(u)h(u)du}+\Big(\frac{1}{h(x)}\Big)\Big(\frac{\int_{0}^{x}w(u)h(u)du}{\int_{0}^{x}w(u)du}\Big)=2.
\end{displaymath}
Hence, $$\Big(h(x)\int_{0}^{x}w(u)du-\int_{0}^{x}w(u)h(u)du\Big)^{2}=0,$$ and this gives $\frac{d}{dx}h(x)=0$ as $\int_{0}^{x}w(u)du\neq 0.$ This completes the proof. $\hfill\Box$\\

Note that $A^{w}(x)=c$ for all $x>0$ characterizes exponential distribution, and so is true for $G^{w}(\cdot)$ and $H^{w}(\cdot).$ If we simultaneously  peep into the lines in the proof of Proposition \ref{pp}, we conclude that $(i),$ $(ii)$ and $(iii)$ get reduced to $A^{w}(x)=G^{w}(x)=H^{w}(x)=c$ for all $x>0.$  \\

In the next proposition we obtain simple relationships between w-AFR, w-GFR and w-HFR functions and hazard rate, that characterize the underlying distributions through their hazard rates. The proof is omitted.
\begin{p1}
\label{propo}
Let $h(x)$ be differentiable for all $x\geq 0.$ Then for any non-negative weight function $w(x),$ and for suitable positive values of constants, $a,b,c, k$ we have
\begin{enumerate}
\item[(i)] $A^{w}(x)=ah(x)$ for all $x$ if and only if $h(x)=k\Big(\int_{0}^{x}w(u)du\Big)^{(1-a)/a}$
\item[(ii)] $G^{w}(x)=bh(x)$ for all $x$ if and only if $h(x)=k\Big(\int_{0}^{x}w(u)du\Big)^{\ln(e/b)-1}$
\item[(iii)] $H^{w}(x)=ch(x)$ for all $x$ if and only if $h(x)=\Big(\frac{1}{kc}\int_{0}^{x}w(u)du\Big)^{1-c},$ where $k$ is an arbitrary constant.
$\hfill\Box$
\end{enumerate}
\end{p1}

One can wonder whether for any particular class of well known probability distribution, weighted means are proportional to their respective hazard rates. If we choose weight function as $w(x)=x^{n},$ then the corresponding failure rate is that of two-parameter Weibull distribution with shape parameter $(n-an+1)/a$ and scale parameter $ka/\big((n-an+1)(n+1)^{\frac{1-a}{a}}\big),$ provided $\Big(\frac{1+n-an}{a}\Big)>0.$  Intuitively, it follows that $(ii)$ and $(iii)$ are also satisfied for two-parameter Weibull distribution having a different sets of scale and shape parameters. We summarize this discussion by claiming that $w(x)=x^{n}$ for suitable $n,$ and $x>0$ is a proper choice of weight function as it results in a legitimate probability distribution. A little work out will show that if we choose $w(x)=e^{nx},$  then the  resultant failure rate function is $h(x)$ does not correspond to a well defined probability distribution, underlying the fact that proportionality of weighted means and hazard rate do not hold good. \\ 
\begin{comment}
{\bf Proof.} We note that $A^{w}(x)=a r(x)$ implies $$\frac{\frac{d}{dx} r(x)}{r(x)}=\Big(\frac{1-a}{a}\Big)\frac{w(x)}{\int_{0}^{x}w(u)du},$$ proving $(i).$ Also, $G^{w}(x)=a r(x)$ implies $$\frac{\frac{d}{dx} r(x)}{r(x)}=(-\ln b)\frac{w(x)}{\int_{0}^{x}w(u)du},$$ proving $(i).$
\end{comment}

We end this subsection by stating some crucial observations in the upcoming remark.
\begin{r1}An essence of introducing the weighted version of means of failure rate lies in the aforementioned Proposition \ref{propo}, where a suitable choice of weight function  characterizes some well known distributions. The readers may also note that, proportionality of each of $A^{w}(\cdot), G^{w}(\cdot)$ and $H^{w}(\cdot)$ with $h(\cdot)$ imply that $h(x)$ is increasing (decreasing) in $x$ if and only if $a\leq (\geq) 1,$  $b\leq (\geq) 1,$ and $c\leq (\geq) 1$  respectively. It is clear that, under the aforementioned conditions, monotonicity of $h(\cdot)$ is independent of the choice of weight function.
\end{r1}
\subsection{Bounds and Limiting behavior of aging means}
${}$\hspace{0.8cm}We state a result from Wijsman (1985) in the form of a lemma. 
\begin{l1}
\label{imple}
Let  $f_i, g_i$  are  non-negative functions, such that the  integrals $\int f_i g_i$ are positive for $i,j=1,2.$ Then
\begin{equation}
\label{newin}
\frac{\int f_1g_1 d\mu }{\int f_1g_2 d\mu}\geq \frac{\int f_2 g_1 d\mu }{\int f_2 g_2 d\mu},\nonumber
\end{equation}
provided $f_{1}/f_{2}$ and $g_1/g_2$ are monotonic in same direction. The inequality in (\ref{newin}) is reversed if $f_{1}/f_{2}$ and $g_1/g_2$ are monotonic in opposite direction. Equality holds if and only if either $f_1/f_2$ or $g_1/g_2$ is a constant. Here $\mu$ is Lebesgue measure on a subset of the real line or counting measure on a subset of the integers.
\end{l1}

The following proposition decides bounds of the aging means on the basis of monotonicity of weight function and hazard rate (as the case may be). An interpretation of the result is 
\begin{p1}\label{che}
\begin{itemize}
    \item[(i)] $A^{w}(x)\geq (\leq)A(x)$ according as $w(x)$ and $h(x)$ are monotonic in same (opposite) direction.
    \item[(ii)] If the hazard rate function $h(x)\geq 1$ for all $x\geq 0$ then $G^{w}(x)\geq (\leq)G(x)$ according as $w(x)$ and $h(x)$ are monotonic in same (opposite) direction.
    \item[(iii)]$H^{w}(x)\geq (\leq)H(x)$ according as $w(x)$ and $h(x)$ are monotonic in same (opposite) direction.
\end{itemize}
\end{p1}
{\bf Proof.} We choose  $f_1(x)=w(x), g_{1}(x)= h(x), f_{2}(x)=g_{2}(x)=1,$ to prove $(i).$ By choosing $f_1(x)=w(x), g_{1}(x)=\ln h(x), f_{2}(x)=g_{2}(x)=1,$ and assuming $\ln h(x)\geq 0$ (since $h(x)\geq 1$ for all $x\geq 0$), Lemma \ref{imple} gives $\Big(\frac{\int_{0}^{x}w(u)\ln h(u)du} {\int_{0}^{x}w(u)du}\Big) \geq (\leq)\Big(\frac{1}{x}\int_{0}^{x}\ln h(u)du\Big)$ according as $w(x)$ and $h(x)$ are monotonic in same (opposite) direction. This proves $(ii).$ Similarly, taking $f_1(x)=w(x), g_{1}(x)= 1, f_{2}(x)=1, g_{2}(x)=1/h(x),$ we prove $(iii).$$\hfill\Box$\\

The readers may arrive at the following remark by looking at the Proposition \ref{che} and the fact that $A^{w}(x)\geq G^{w}(x) \geq H^{w}(x)$ for all $x\geq 0.$
\begin{r1}
If $w(x)$ and $h(x)$ are monotonic in same direction then the lower and upper bounds of the aging means of failure rate are $H(x)$ and $A^{w}(x)$ respectively. On the other hand, $w(x)$ and $h(x)$ are monotonic in opposite direction then the lower and upper bounds of the aging means of failure rate are $H^{w}(x)$ and $A(x)$ respectively. The lower and upper bounds of the aging means of failure rate discussed in this article are $\min(H(x),H^{w}(x))$ and $\max(A(x),A^{w}(x))$ respectively. 
\end{r1}

In the following theorem, we obtain bounds for the ratio of the weighted hazard means by associating weights in sequence.

\begin{t1}
\label{inter}
Let $h_{k}(x)=w(x)h_{k-1}(x)=(w(x))^{k}h(x), k\geq 1, h_{0}(x)=h(x),$ for $x>0.$ We define $A^{w}_{h_{k}}(x)=\Big(\frac{\int_{0}^{x}w(u)h_{k}(u)du}{\int_{0}^{x}w(u)du}\Big),$ $G^{w}_{h_{k}}(x)=\exp\Big(\frac{\int_{0}^{x}w(u)\ln h_{k}(u)du}{\int_{0}^{x}w(u)du}\Big),$ and $H^{w}_{h_{k}}(x)=\Big(\frac{\int_{0}^{x}w(u)du}{\int_{0}^{x}\frac{w(u)}{h_{k}(u)}du}\Big).$ For $x>0,$ 
 the following statements hold.
\begin{enumerate}
\item[(i)] If $h(x)$ and $w(x)$ are monotonic in opposite (same) direction then $$\frac{A^{w}_{h_{k}}(x)}{A^{w}_{h}(x)}\geq (\leq) \frac{\int_{0}^{x}w(u)du}{\int_{0}^{x}(w(u))^{n+1}du}.$$
\item[(ii)] If $w(x)>1,$ then $$\frac{G^{w}_{h_{k}}(x)}{G^{w}_{h}(x)}\geq \exp\Big(\frac{k}{x}\int_{0}^{x}\ln w(u)du\Big).$$
\item[(iii)] If $h(x)$ and $w(x)$ are monotonic in same (opposite) direction then $$\frac{H^{w}_{h_{k}}(x)}{H^{w}_{h}(x)}\geq (\leq)\frac{\int_{0}^{x}w(u)du}{\int_{0}^{x}\frac{1}{(w(u))^{k-1}}du}.$$
\item[(iv)] If $w(x)$ and $h(x)$ are monotonic in same (opposite) direction then $A^{w}_{h_{k}}(x)\geq (\leq) A(x),$ and $H^{w}_{h_{k}}(x)\geq (\leq) H(x),$ according as $w(x)\leq (\geq) 1.$
\item[(v)]If $h(x)\geq 1,$ $w(x)\geq 1$  then $G^{w}_{h_{k}}(x)\geq G(x),$ provided $w(x)$ and $h(x)$ are monotonic in same direction. 
\end{enumerate}
\end{t1}
{\bf Proof.} The proofs of $(i), (ii)$ and $(iii)$ follow by applying Lemma \ref{imple} on the ratios, viz., $$\frac{A^{w}_{h_{k}}(x)}{A^{w}_{h}(x)}=\frac{\int_{0}^{x}(w(u))^{k+1}h(u)du}{\int_{0}^{x}w(u)h(u)du},\frac{G^{w}_{h_{k}}(x)}{G^{w}_{h}(x)}=\exp\Big(\frac{k\int_{0}^{x}w(u)\ln w(u)du}{\int_{0}^{x}w(u)du}\Big),$$ and  $$\frac{H^{w}_{h_{k}}(x)}{H^{w}_{h}(x)}=\Big(\frac{\int_{0}^{x}\frac{w(u)}{h(u)}du}{\int_{0}^{x}\frac{1}{w^{k-1}(u)h(u)}du}\Big), x>0.$$ The proof of $(iv)$ follows from $(i)$ and $(iii)$ of Proposition \ref{che}. The proof of $(v)$ follows from $(ii)$ of Proposition \ref{che}.     If $w(x)$ and $h(x)$ are monotonic in same (opposite) direction then $A^{w}_{h_{k}}(x)\geq (\leq) A^{w}_{h}(x),$ and $H^{w}_{h_{k}}(x)\geq (\leq) H^{w}_{h}(x),$ according as $w(x)\geq (\leq) 1.$$\hfill\Box$

The above theorem can be interpreted by saying that one can keep minimizing the means of failure rate (AFR and HFR) of a component, having increasing failure rate by associating weights in sequence which are monotonically decreasing with time. However, GFR increases rapidly with the increase in number of weight functions and is independent of the nature of monotonicity of weight and hazard rate. Theorem \ref{inter} $(ii)$ reveals that if $k\rightarrow \infty$ and $w(x)>1$ then $G^{w}_{h_{k}}(x)\rightarrow \infty.$

\section{Non-parametric classes of distributions based on weighted means of failure rates}

${}$\hspace{0.8cm}We define non-parametric classes of distributions on the basis of monotonicity of w-AFR, w-GFR and w-HFR. 
\begin{d1}
\label{define}
A random variable $X$ is said to belong to the class of
\begin{enumerate}
 \item[$(i)$] Increasing (resp. Decreasing) weighted arithmetic mean failure rate $Iw-AFR$ (resp. $Dw-AFR$)) distributions if  $A^{w}(x)$ is increasing (resp. decreasing) in $x>0.$
\item[$(ii)$] Increasing  (resp. Decreasing) weighted geometric mean failure rate $Iw-GFR$ (resp. ($Dw-GFR$)) distributions if $G^{w}(x)$ is increasing (resp. decreasing) in $x>0.$
\item[$(iii)$] Increasing (resp. Decreasing)  weighted harmonic mean failure rate  $Iw-HFR$ (resp. $Dw-HFR$)) distributions if  $H^{w}(x)$
 is increasing (resp. decreasing) in $x>0.$\hfill$\Box$
 \end{enumerate}
 \end{d1}
 ${}$\hspace{0.8cm}The next theorem emphasises on the fact that the monotonic behaviour of $h(x)$ is possessed by $A^{w}(x), G^{w}(x)$ and $H^{w}(x).$
\begin{t1}
\label{incr}
If $h(x)$ is increasing (decreasing) in $x\geq 0$  then 
\begin{enumerate}
\item[(i)]$A^{w}(x)$ is increasing (decreasing) in $x\geq 0;$ 
\item[(ii)]$G^{w}(x)$ is increasing (decreasing) in $x \geq 0;$ 
\item[(iii)] $H^{w}(x)$ is increasing (decreasing) in $x\geq 0.$ 
\end{enumerate}
\end{t1}
{\bf Proof.} To prove $(i)$, we note that $\Big(\int_{0}^{x}w(u)du\Big)\Big(\frac{d}{du}A^{w}(x)\Big)=w(x)(h(x)-A^{w}(x)),$ 
and thus $\Big(\frac{d}{dx}A^{w}(x)\Big) \geq (\leq)~ 0$ according as $h(x) \geq (\leq) A^{w}(x)$ for all $x\geq 0.$  If $h(x)$ is increasing (decreasing) in $x$ then $h(x)\geq (\leq)~ A^{w}(x)$ for $x\geq 0.$  This proves $(i).$ 
Similarly, to prove $(ii),$ we first note that $$\Big(\frac{d}{dx}G^w(x)\Big)=\frac{G^w(x)}{\Big(\int_{0}^{x}w(u)du\Big)}w(x)\ln \Big(\frac{h(x)}{G^{w}(x)}\Big),$$ and this implies that $\frac{d}{dx}G^{w}(x)\geq (\leq) ~0$ according as $h(x)\geq(\leq)G^{w}(x).$ One can note that if $h(x)$ is increasing (decreasing) in $x$ then $h(x)\geq (\leq)~G^{w}(x)$ for all $x\geq 0,$ thus proving $(ii).$ To prove $(iii),$ we first note that 
$$\Big(\frac{d}{dx}H^{w}(x)\Big)\Big(\int_{0}^{x}\frac{w(p)}{h(p)}dp\Big)=w(x)\Big\{1-\frac{H^{w}(x)}{h(x)}\Big\},$$ 
and hence we find that $\Big(\frac{d}{dx}H^{w}(x)\Big) \geq (\leq)~ 0$ according as $h(x) \geq (\leq) H^{w}(x)$ for all $x\geq 0.$ Also, if $h(x)$ is increasing (decreasing) in $x$ then $h(x)\geq (\leq) H^{w}(x)$ for all $x\geq 0.$ This completes the proof.$\hfill\Box$\\

The next example highlights the importance of choosing  weight functions in generating new distributions. It is also cited in upcoming counterexample \ref{impcou} for establishing that w-IFRA class is not closed under formation of coherent systems.
 \begin{ex}
\label{impex}
Let $X$ follows two parameter Weibull distribution with scale and shape parameter $\alpha$ and $\beta$ respectively. If we take $w(x)=e^{nx}$ for all $x>0,$ then the weighted random variable $X^{w}$ has failure rate $h^{w}(x)=\alpha \beta e^{nx} x^{\beta-1}.$  Here, taking $n=-m$ with $m>0$,
\begin{eqnarray}
\int_{0}^{x}w(u)h(u)u
&=&\alpha \beta\int_{0}^{x}e^{-mt} t^{\beta-1}dt\nonumber\\
&=&\alpha\beta(-n)^{-\beta} \gamma(\beta,-nx)\nonumber\\
\label{who1}
&=&\alpha\beta(-n)^{-\beta} \Big(\Gamma(\beta)-\Gamma(\beta,-nx)\Big),
\end{eqnarray}
where the incomplete Gamma function $\gamma (z,a)$ and its complement $\Gamma(z,\alpha)$ (also known as Prym's function) are 
\begin{equation}
\gamma(a,x)=\int_{0}^{x} t^{a-1}e^{-t}dt, \Gamma(a,x)=\int_{x}^{\infty} t^{a-1}e^{-t}dt,~ \mbox{Real} (a)>0),\nonumber
\end{equation}
satisfying $\gamma(a,x)+\Gamma(a,x)=\Gamma(a).$ If $n<0$ we have real values for $\bar{F}^{w}(t),$ as
$$\bar{F}^{w}(x)=\exp\Big\{-\alpha\beta(-n)^{-\beta} \Big(\Gamma(\beta)-\Gamma(\beta,-nx)\Big)\Big\}, x>0, \beta>0.$$ 
Also, considering $m=-n,$ we get
\begin{equation}
\frac{d}{dx} h^{w}(x)=\frac{d}{dx}\Big(\alpha\beta e^{nx}x^{\beta-1}\Big)=\alpha \beta e^{nx}x^{\beta-2}(\beta-1-mx) \leq 0\nonumber
\end{equation}
if $(\beta-1-mx)\leq 0,$ i.e., $\frac{d}{dx} h^{w}(x)\leq 0$ if $x\geq \frac{\beta-1}{m}.$ If $\beta<1$ then $\frac{d}{dx} h^{w}(x)\leq 0$ for all $x>0.$ Thus, $X^{w}$ is DFR if $\beta<1.$
On the other hand if $\beta>1,$ then $\frac{d}{dx} h^{w}(x)\geq  0$ for $x \in (0,\frac{\beta-1}{m})$ and $\frac{d}{dx} h^{w}(x)\leq 0$ for $x \geq \frac{\beta-1}{n}.$
Thus $X^{w}$ is DFR if $\beta<1,$ whereas $X^{w}$ has upside-down bathtub shaped failure rate if $\beta>1.$ 
Using (\ref{who1}) and the fact that $\int_{0}^{x}w(u)du=\frac{1}{n}\Big(e^{nx}-1\Big)$ we get
\begin{equation}A^{w}(x)
\label{who2}
=\frac{n (-n)^{-\beta}\alpha \beta \Big(\Gamma[\beta]-\Gamma[\beta,-nx]\Big)}{\Big(e^{nx}-1\Big)}\nonumber
\end{equation}
Here, for $\beta<1,$ $h^{w}(x)$ is decreasing in $x,$ and so is $A^{w}(x)$ as evident from Theorem \ref{incr}. Similarly, for $\beta>1,$ $A^{w}(x)$ is upside-upside-down bathtub.
Using (\ref{who2}), we note that, 
\begin{eqnarray}
\label{lq}
\frac{d}{dx} A^{w}(x)&=&n(-n)^{-\beta}\alpha\beta \frac{d}{dx}\Big(\frac{\Big(\Gamma[\beta]-\Gamma[\beta,-nx]\Big)}{e^{nx}-1}\Big)\nonumber\\
&=& n(-n)^{-\beta}\alpha\beta \Bigg\{\frac{(e^{nx}-1)e^{nx}(-nx)^{\beta-1}(-n)-\Big(\Gamma[\beta]-\Gamma[\beta,-nx]\Big)e^{nx}n}{(e^{nx}-1)^{2}}\Bigg\}\nonumber\\
&=&n(-n)^{-\beta}\alpha\beta \Bigg\{\frac{(e^{nx}-1)e^{nx}(-nx)^{\beta}\Big(\frac{(-n)}{-(nx)}\Big)-\Big(\Gamma[\beta]-\Gamma[\beta,-nx]\Big)e^{nx}n}{(e^{nx}-1)^{2}}\Bigg\}\nonumber\\
&=&\frac{n(-n)^{-\beta}\alpha\beta e^{nx}}{x(e^{nx}-1)^{2}}\Big\{(e^{nx}-1)(-nx)^{\beta}+nx\Big(\Gamma[\beta,-nx]-\Gamma[\beta]\Big)\Big\}\nonumber\\
%&=&\frac{1}{(e^{nt}-1)^{2}}e^{nt}nt^{\beta-1}(-nt)^{-\beta}\alpha \beta \Big\{(e^{nt}-1)(-nt)^{\beta}+ nt\Big(-\Gamma[\beta]+\Gamma[\beta,-nt]\Big)\Big\}\nonumber\\
&=&\frac{n^{2}x(-n)^{-\beta}\alpha\beta e^{nx}}{x(e^{nx}-1)^{2}}\Big\{\Gamma[\beta,-nx]-\Gamma[\beta]-(e^{nx}-1)(-nx)^{\beta-1}\Big\}\nonumber\\
&=&\frac{n^{2}x(-n)^{-\beta}\alpha\beta e^{nx}}{x(e^{nx}-1)^{2}}\Big\{-\gamma[\beta,-nx]-P(x)\Big\},\mbox{(say)},\nonumber
\end{eqnarray}
where $$P(x)=(e^{nx}-1)(-nx)^{\beta-1}\leq 0, x>0.$$
The change point of monotonicity of $A^{w}(x)$ is determined by the root of equation $\gamma[\beta,-nx]+P(x)=0.$
Similarly, we obtain
$$G^{w}(x)=\alpha  \beta  t^{\beta -1} (-n t)^{\frac{\beta -1}{e^{n t}-1}} e^{\frac{(\beta -1) (E_1(-n t)+\gamma )}{e^{n t}-1}}$$
$$H^{w}(x)=\frac{\alpha \beta (e^{nt}-1)nt^{\beta}(-nt)^{-\beta}}{\Gamma[2-\beta]-\Gamma[2-\beta,-nt]}, $$ where $\gamma \sim 0.577216$ is Euler's constant and  $E_{n}(z)$ is the exponential integral function.
\end{ex}
${}$\hspace{0.8cm} Below, we state two theorems highlighting the inclusion property of the non-parametric aging classes given in Definition \ref{define}. The proof follows due to Theorem \ref{incr}, line of the proof therein and the fact that $A^{w}(x)\geq G^{w}(x) \geq H^{w}(x)$ for all $x>0.$ 
\begin{t1}\label{theo3.1}
	$IFR\subseteq Iw-AFR\subseteq Iw-GFR\subseteq Iw-HFR$
\end{t1}
\begin{t1}\label{theo3.2}
	$DFR\subseteq Dw-HFR\subseteq Dw-GFR\subseteq Dw-AFR$
\end{t1}
\subsection{Characterization results for w-AFR, and w-GFR}
${}$\hspace{0.8cm} We introduce the concept of weighted star-shaped (anti-star) function which is a generalization of star-shaped (anti-star) function to give an equivalent condition of  $Iw-AFR$ and $Dw-AFR$ classes of distributions.
\begin{d1}
A function $g(x)$ defined on $[0,\infty)$ is said to be weighted star-shaped (weighted anti-star shaped) with respect to a non-negative weight function $w(x)$ if $\Big(-\frac{1}{\int_{0}^{x}w(u)du}\Big)~g(x)$ is increasing in $x>0.$ Equivalently, for $0\leq \alpha \leq 1$ and $x \geq 0,$ 
$$g(\alpha x) \leq \Big(\frac{\int_{0}^{\alpha x}w(u)du}{\int_{0}^{x}w(u)du}\Big) g(x)$$
\end{d1}
${}$\hspace{0.8cm}The next theorem gives a necessary and sufficient condition of a increasing (decreasing) weighted arithmetic failure rate or weighted failure rate average class of distributions, denoted by w-AFR. We omit the proof for the sake of brevity.
\begin{t1}
Let $X$ has increasing (decreasing) w-AFR. 
 Then the following conditions are equivalent.
\begin{enumerate}
    \item[(i)] $\Big(-\frac{1}{\int_{0}^{x}w(u)du}\Big)\ln \bar{F}^{w}(x)$ is increasing (decreasing) in $x>0.$ 
    \item[(ii)] $-\ln \bar{F}^{w}(x)$ is weighted star-shaped  (weighted anti-star shaped) with respect to $c(\cdot)$.
    \item[(iii)] $\Big(\bar{F}^{w}(x)\Big)^{\frac{1}{\int_{0}^{x}w(u)du}}$ is decreasing (increasing) in $x>0.$
    \item[(iv)] For $\alpha \in [0,1],$ and $x>0,$ $\bar{F}^{w}(\alpha x)\geq (\leq) \Big(\bar{F}^{w}(x)\Big)^{\frac{\int_{0}^{\alpha x}w(u)du}{\int_{0}^{x}w(u)du}}.$ \hfill$\Box$
\end{enumerate}
\end{t1}
${}$\hspace{0.8cm}The following theorem gives equivalent conditions for $w-GFR$. 
\begin{t1}
Let $X$ has increasing (decreasing) w-GFR.Then the following conditions are equivalent.
\begin{enumerate}
    \item[(i)] $\Big(\frac{1}{\int_{0}^{x}w(u)du}\Big)\Big(\int_{0}^{x}w(u)\ln h(u)du\Big)$ is increasing (decreasing) in $x>0.$ 
    \item[(ii)] $\Big(\int_{0}^{x}w(u)\ln h(u)du\Big)$ is weighted star-shaped  (weighted anti-star shaped) with respect to $c(\cdot)$.
       \item[(iii)] For $\alpha \in [0,1],$ and $x>0,$ $\int_{0}^{\alpha x}w(u)\ln h(u)du\leq (\geq) \frac{\int_{0}^{\alpha x}w(u)du}{\int_{0}^{x}w(u)du}\Big(\int_{0}^{x}w(u) \ln h(u)du\Big).$ \hfill$\Box$
\end{enumerate}
\end{t1}
We note that, $0\leq \frac{\int_{0}^{\alpha x}w(u)du}{\int_{0}^{x}w(u)du}\leq 1$ since $w(x)\geq 0$ for all $x\geq 0$ and $0\leq \alpha \leq 1.$
\subsection{Results on coherent system}
In this section, we primarily focus on Iw-AFR class and its closure properties. 
We know that IFRA class is closed under the formation of coherent system. Naturally, a question ponders, whether the same result is true for increasing weighted AFR class (Iw-AFR). Let us consider a coherent system with $n$ components having weighted survival functions $\bar{F}_{i}^{w}(x)$ for $i=1,2,\ldots, n.$ The survival function $\bar{F}^{w}(x)$ of the resultant coherent system satisfies \begin{equation}
\label{coh}
\bar{F}^{w}(\alpha x)=h(\bar{F}_{1}^{w}(\alpha x), \bar{F}_{2}^{w}(\alpha x), \ldots, \bar{F}_{n}^{w}(\alpha x)),
\end{equation}
where $h$ represents the survival function of the coherent system.
Further, if we assume that each $X_{i}$ has increasing w-AFR, then we explore what would be the survival function of the resultant coherent system. Since  $\bar{F}_{i}^{w}(\alpha x)\geq \Big(\bar{F}_{i}^{w}(x)\Big)^{\frac{\int_{0}^{\alpha x}w(u)du}{\int_{0}^{x}w(u)du}}$ for $i=1,2,\ldots,n,$ $\alpha \in [0,1], x>0,$ and  $h$ is increasing in each argument, (\ref{coh}) reduces to 
\begin{equation}
\label{coh1}
\bar{F}^{w}(\alpha t)\geq h\Big(\big(\bar{F}_{1}^{w}(x)\big)^{\frac{\int_{0}^{\alpha x}w(u)du}{\int_{0}^{x}w(u)du}}, \big(\bar{F}_{2}^{w}(x)\big)^{\frac{\int_{0}^{\alpha x}w(u)du}{\int_{0}^{x}w(u)du}}, \ldots, \big(\bar{F}_{n}^{w}(x)\big)^{\frac{\int_{0}^{\alpha x}w(u)du}{\int_{0}^{x}w(u)du}}\Big)\nonumber
\end{equation}
The following counterexample shows that Iw-AFR is not closed under the formation of coherent system.
\begin{co}
\label{impcou}
Let us consider a series system with lifetime $X$ formed by two components with lifetimes $X^{w}_{1}$ and $X^{w}_{2}$ respectively. Let the  failure rates be $h_{1}(t),$ and $h_{2}(t)$  with corresponding weights $w_{1}(t)$ and $w_{2}(t)$ respectively. Let $h_{1}(t)=\alpha\beta t^{\beta-1}, w_{1}(t)=e^{nt}, $ and $h_{2}(t)=abt^{b-1}, w_{2}(t)=(1-e^{nt})$ where $\alpha, a>0; \beta, b >1; n<0.$
Since, $\beta, b >1;$ $h_{1}(t)$ and $h_{2}(t)$ are increasing in $t.$ By Theorem \ref{incr}, it follows that
$A^{w}_{1}(t)$ and $A^{w}_{2}(t)$ are increasing in $t$ as $h_{1}(t)$ and $h_{2}(t)$ are increasing in $t.$  Then the hazard rate of the series system is given by $h_{X}(t)=h_{1}(t)w_{1}(t)+h_{2}(t)w_{2}(t)$ for all $t>0.$ From Example \ref{impex}, it follows that each of $h_{1}^{w}(t)=h_{1}(t)w_{1}(t)$ and $h_{2}^{w}(t)=h_{2}(t)w_{2}(t)$ are non-monotonic in $t>0.$ (upside-down bathtub curve). Thus, $X^{w}_{1}$ and $X^{w}_{2}$ are Iw-AFR but not IFR. Here, $X$ is not Iw-AFR since $h(t)$ is non-monotonic (as noted in Example \ref{impex}) and non-monotonicity of $h(t)$ is transmitted to $A(t)$ (by Theorem \ref{incr}).\end{co}
\begin{comment}
as evident from 
\begin{eqnarray}
\frac{d}{dt}A(t)&=& \frac{d}{dt}\Bigg\{\frac{1}{t}\int_{0}^{t}\Big(\alpha\beta u^{\beta-1}e^{nu}+a\beta u^{b-1}(1-e^{nu})\Big)du\Bigg\}\nonumber\\
&=&\frac{-a t^b \left(b \left(e^{n t}-1\right)+1\right)+a b t^b (-n t)^{-b} (\Gamma (b)-\Gamma (b,-n t))+\alpha  \beta  e^{n t} t^{\beta }-\alpha  \beta  t^{\beta } (-n t)^{-\beta } (\Gamma (\beta )-\Gamma (\beta ,-n t))}{t^2}.\nonumber
\end{eqnarray} 
In particular,  if we tacitly choose the values of the parameters, $\frac{d}{dt}A(t)$ changes sign, making $A(t)$ a non-monotonic function. Notably, if choose $b<\beta,$
\end{comment}
\begin{comment}
\begin{l1}
Since $f(x)=x^{\psi(x)}$ is concave in $x$ if $0\leq \psi(x)\leq 1$ for all $x\geq 0.$
\end{l1}
\end{comment}
\begin{comment}
\begin{p1}
Let $X$ and $X^{w}$ be two random variables having failure rates $r_{X}(t)$ and $r_{X}^{w}(t)=w(t)r_{X}(t)$ respectively for all $t>0.$ Let the weighted means of $X$ be denoted by $A_{X}^{w}(t),G_{X}^{w}(t),$  and $H_{X}^{w}(t),$ respectively. On the other hand, weighted means of $Y$ be denoted by $A_{Y}(t),G_{Y}(t),$  and $H_{Y}(t),$ respectively.  Then
$A_{X}^{w}(t)=A_{Y}^{w}(t),$ $G_{X}^{w}(t)=G_{Y}(t)\exp\Big(-\frac{\int_{0}^{x}w(u)r(u)du}{\int_{0}^{x}w(u)du}\Big),$
\end{p1}   
\end{comment}
${}$\hspace{0.8cm}We end this section by the following remark which once again emphasizes the importance of the concepts introduced in this article which in turn engenders the upcoming section.
\begin{r1}
\label{rem1}
\begin{comment}
The monotonic behaviour of $w(\cdot)$ is transmitted to $h(\cdot)$ in case each of w-AFR, w-GFR and w-HFR is proportional to hazard rate.
\end{comment}
 In particular, if in the definition \ref{newd} discussed above, we replace failure rate $h(\cdot)$ by hazard quantile function $h_q(\cdot)$ and $w(\cdot)$ by density quantile function $q(\cdot)$  respectively, with support restricted to $[0,1],$ we get  quantile version of AFR, GFR and HFR as mentioned in upcoming section.
\end{r1}
 \section{Quantile version of AFR, GFR and HFR}
${}$\hspace{0.8cm}The Remark \ref{rem1} is the genesis of this section. The readers would be interested to see how the act of replacing failure rate by  hazard quantile function and weight function by density quantile function respectively will differ the rest of the analysis and is taken up in the present section.
${}$\hspace{0.8cm}For a random variable $X$, quantile function (QF) is defined as 
\begin{equation}
    Q(u)=F^{-1}(u)=\inf\Big\{x: F(x)\geq u\Big\}, 0\leq u\leq 1
\end{equation}
gives $FQ(u)=u.$ Differentiating with respect to $u,$ we get $f(Q(u))q(u)=1$ or $f(Q(u))=\frac{1}{q(u)}$, where $f(Q(u))$ and $q(u) = \frac{d}{du}Q(u)$ are respectively known as the density quantile function and quantile density function of $X$.  From the definition of hazard rate, the corresponding hazard quantile function is given by
\[h_q(u) = h(Q(u)) = \frac{f(Q(u))}{\bar{F}(Q(u))} = \frac{1}{(1-u)q(u)}.\] 
This implies $q(u)=\frac{1}{(1-u)h_q(u)}.$ Integrating, we get $Q(u)=\int_{0}^{u}\frac{1}{(1-p)h_q(p)}dp $.  The quantile approach is an alternative to the traditional distribution function method as it can also used to specify a probability distribution.  As the quantile approach possess some interesting properties not shared by its distribution function counterpart and in many situations, quantile measures provide simple expressions that are easily amenable to many computational analysis.  Abundant literature are now available on various properties of quantile functions and different measures based on it and their applications, for details see Gilchrist (2000), Nair et al. (2013), Nair et al. (2023) and Aswin et al. (2023) and references therein.  \\

Note that, the quantile version of AFR, GFR and HFR, denoted by $QA(\cdot), QG(\cdot), QH(\cdot)$ respectively, can be independently derived using quantile approach. Alternatively, we can obtain the same by applying the Remark \ref{rem1} (ii) in Definition \ref{newd}, $i.e.$, replacing the failure rate $h(\cdot)$ by hazard quantile function $h_q(\cdot)$ and $w(\cdot)$ by density quantile function $q(\cdot)$  respectively, with support restricted to $[0,1].$  
\begin{eqnarray}\label{2.2}
 QA(u)=QA(Q(u))&=&\frac{-\ln (1-F(Q(u)))}{Q(u)} \nonumber\\
 &=& \frac{-\ln (1-u)}{Q(u)} = -\Big(\ln (1-u)\Big)\Big\{\int_{0}^{u}\frac{1}{(1-p)h_q(p)}dp\Big\}^{-1}.\label{afr}
  \end{eqnarray}
 \begin{eqnarray}
QG(u)=QG(Q(u))&=&\exp\Big(\frac{1}{Q(u)}\int_{0}^{u}\ln\Big(\frac{1}{(1-p)q(p)}\Big) dQ(p)\Big)\nonumber\\ &=&\exp\Big(-\frac{1}{Q(u)}\int_{0}^{u}q(p)\ln\Big((1-p)q(p)dp\Big)\Big),\nonumber
  \end{eqnarray}
 or equivalently,
 \begin{eqnarray}
 \label{impg}
 Q(u)\ln QG(u)&=&-\int_{0}^{u}(\ln (1-p))q(p)dp-\int_{0}^{u}(\ln q(p))q(p)dp\nonumber\\
\label{gfr1} &=&-\int_{0}^{u}q(p)\ln \Big\{(1-p)q(p)\Big\}dp=\int_{0}^{u}q(p)\ln h_q(p) dp
  \end{eqnarray}
%\end{comment}
and
\begin{eqnarray}\label{2.4}
QH(u)=QH(Q(u))&=&\Big(\frac{1}{Q(u)}\int_{0}^{u}\frac{1}{h(Q(p))}dQ(p)\Big)^{-1}\nonumber\\
&=&Q(u)\Big(\int_{0}^{u}(1-p)(q(p))^{2}dp\Big)^{-1}\\
&=&Q(u)\Big(\int_{0}^{u}\frac{q(p)}{h_q(p)}dp\Big)^{-1}.
\end{eqnarray}
Differentiating \eqref{2.2} with respect to $u$, we obtain
\[
QA'(u) Q(u) + QA (u) q(u) = \frac{1}{1 - u}
\]
When quantile AFR is increasing (decreasing), we get 
\[
QA(u) \leq (\geq) ~h_q(u).
\]
From \eqref{2.4}, we have
\[
\frac{Q(u)}{QH(u)} = \int_{0}^{u} (1 - p) (q(p))^2 dp.
\]
Differentiating with respect to $u$, we get
\[
QH(u) q(u) - Q(u) QH' (u) = (1 - u) (q(u))^2 \left(QH (u) \right)^2.
\]
Now when the quantile HFR is increasing (decreasing), yield
\[
QH(u) \leq (\geq) h_q(u).
\]
A similar arguement as given in Theorem \ref{incr} depicts that monotonicity of hazard quantile function $h_q(\cdot)$ is transmitted to quantile version of AFR, GFR and HFR, $i.e.$, $QA(\cdot), QG(\cdot)$ and $QH(\cdot).$ In continuation to the Proposition \ref{propo} of previous section, if we replace $w(\cdot)$ by $q(\cdot)$ and $h(\cdot)$ by $h_q(\cdot)$, we find that proportionality of weighted means of quantile hazard functions with quantile hazard function characterizes some quantile function. To the best of our knowledge,  $Q(x)$ as obtained in Proposition \ref{propo1} represents a new generalized version of quantile function where $Q(0)\neq 0.$\\

The next example gives the quantile function of AFR, GFR and HFR of Pareto-I distribution.
\begin{ex}
For Pareto I distribution, with quantile function $Q(u)=\alpha (1-u)^{-1/\alpha},$ we have $$QA(u)=-\frac{(1-u)^{1/\alpha } \log (1-u)}{\alpha }, QG(u)=e^{1-(1-u)^{1/\alpha }} (1-u)^{1/\alpha }, \mbox{~and ~} QH(u)=-\frac{2 (1-u)^{1/\alpha }}{(1-u)^{2/\alpha }-1},$$ for $0<u<1.$
\end{ex}
\begin{p1}
\label{propo1}
Let hazard quantile function $h_q(u)$ be differentiable for all $u\in [0,1].$ Then for the non-negative weight function $q(u),$  called as density quantile function and for $a,b,c>0$ we have
\begin{enumerate}
\item[(i)] $QA(u)=a \; h_q(u)$ for all $u\in [0,1]$ if and only if $Q(u)=\Big(\frac{1}{ak}\Big)^{a}\Big\{\ln(\frac{A}{1-u})\Big\}^{a}.$
\item[(ii)] $QG(u)=b \; h_q(u)$ for all $u$ if and only if $Q(u)=\Big(\frac{\ln(e/b)}{k}\Big)^{\frac{1}{\ln(e/b)}}\Big\{\ln(\frac{A}{1-u})\Big\}^{\frac{1}{\ln(e/b)}}.$
\item[(iii)] $QH(u)=c \; h_q(u)$ for all $u$ if and only if $Q(x)=\Big(\frac{\ln(e/b)}{k}\Big)^{\frac{1}{\ln(e/b)}}\Big\{\ln(\frac{A}{1-u})\Big\}^{\frac{1}{\ln(e/b)}}$ where $k$ is an arbitrary constant.
\end{enumerate}
\end{p1}
{\bf Proof.} From Theorem \ref{propo} $(i)$, it follows that $QA(u)=a \; h_q(u)$ is equivalent to $$h_q(u)= k \Big(\int_{0}^{u}q(p)dp\Big)^{(1-a)/a}$$ and since $h_q(u)=\frac{1}{(1-u)q(u)},$ we prove $(i).$ Proofs of $(ii)$ and $(iii)$ are similar. $\hfill\Box$\\

Transformation on a random variable is generally employed to find the best model for a given set of observations.  A simple alternative method to this is to keep the original data as it is and transform the QF to find the best model, using the following property of quantile functions which is not
shared by the distribution function. If $T_X(x)$ is a continuous non-decreasing function then $T_X \left(Q_X (u) \right)$ is the QF of $T_X(X)$ or in symbols
\[
Q_{T(X)}(u) = T \left(Q_X (u)\right).
\]
\begin{t1}
Let $T(\cdot)$ be a continuous non-decreasing and invertible transformation.  Then the quantile versions of AFR, GFR and HFR takes the form
\begin{itemize}
    \item [(i)] $QA_{T(X)}(u) = \frac{- \log (1 - u)}{T \left(Q_X (u)\right)}$,
    \item[(ii)] $QG_{T(X)}(u) = \exp \Big (-\frac{1}{T \left(Q_X (u)\right)} \int_{0}^{u} T' \left(Q_X (p)\right)q(p) \left[\log (1 - p) T' \left(Q_X (p)\right)q(p)\right] dp \Big)$, and
    \item[(iii)] $QH_{T(X)} (u) = T \left(Q_X (u)\right) \Big(\int_{0}^{u} (1 - p) \left(T' \left(Q_X (p)\right)q(p)\right)^2 dp\Big)^{-1}.$ 
\end{itemize}
\end{t1}
\begin{t1}
The following statements are equivalent: $(i)$ $X$ follows Exponential distribution with shape parameter $c$, $(ii)$ $QA(u)=c$, $(iii)$ $QG(u)=c$, and $(iv)$ $QH(u)=c.$$\hfill\Box$
\end{t1}
\begin{r1}
The quantile version is not always equivalent to its distribution function approach.
\end{r1}
\begin{t1}
$QA(u) = \left(Q (u) \right)^{C - 1}$, where $C > 0$ holds if and only $X$ follows Weibull distribution with quantile function $Q(u) = \left(-\log (1 - u)\right)^{\frac{1}{\lambda}}, \; 0 < u < 1, \; \lambda > 0$.$\hfill\Box$
\end{t1}
%\subsection{Quantile-based AFR, GFR and HFR for survival models: Proportional hazards model}

For many models, the distribution function and quantile approches yield similar properties as we have seen in Theorem 4.2, while for certain other cases, it gives different results.  For example, when $X$ and $Y$ satisfy proportional hazard rate model (PHM), we have $h_{q_Y}(u)=\theta h_{q_X}(u),$ or equivalently, we have $\bar{F}_{Y}(u)=\big(\bar{F}_{X}(u)\big)^{\theta}.$ We look at the corresponding AFR, GFR and HFR of $Y$.  It is easy to note that $A_{Y}(x)=\theta~ A_{X}(x)$, $G_{Y}(x)=\theta ~G_{X}(x)$ and $H_{Y}(x)=\theta~ H_{X}(x)$.  To obtain the quantile version of AFR, GFR and HFR under PHM, it is easy to obtain the QF of $Y$, as
\begin{eqnarray}
Q_{Y}(u)=\inf\Big\{x:F_{Y}(x)\geq u\Big\}
%&=&\inf\Big\{x:1-\bar{F}_{Y}(x)\geq u\Big\}\nonumber\\
%&=&\inf\Big\{x:1-u\geq (\bar{F}_{X}(x))^{\theta}\Big\}\nonumber\\
%&=&\inf\Big\{x:F_{X}(x)\geq 1-(1-u)^{1/\theta}\Big\}\nonumber\\
=Q_{X}(1-(1-u)^{1/\theta}),\nonumber
\end{eqnarray}
which in turn obtain the quantile version of AFR for PHM as
\begin{eqnarray}
QA_{Y}(u)=-\frac{\ln (1-u)}{Q_{Y}(u)} =\frac{-\ln (1-u)}{Q_{X}(1-(1-u)^{1/\theta})} = \theta QA_{X}(1-(1-u)^{1/\theta}) \neq&\theta QA_{X}(u),\nonumber
\end{eqnarray}
since \begin{eqnarray}
QA_{X}(1-(1-u)^{1/\theta}) = \frac{-\ln \Big\{\Big\{1-(1-(1-u)^{1/\theta}\Big\}\Big\}}{Q_{X}(1-(1-u)^{1/\theta})} = -\frac{1}{\theta}\frac{\ln (1-u)}{Q_{X}(1-(1-u)^{1/\theta})}.
\end{eqnarray}
%Hence, \begin{eqnarray}
%QA_{Y}(u)=\frac{-\ln (1-u)}{Q_{X}(1-(1-u)^{1/\theta})}&=&\theta\Big(-\frac{1}{\theta}\frac{\ln (1-u)}{Q_{X}(1-(1-u)^{1/\theta})}\Big)\nonumber\\
%&=&\theta QA_{X}(1-(1-u)^{1/\theta})\nonumber\\
%&\neq&\theta QA_{X}(u)
%\end{eqnarray}
The quantile GFR of PHM will be
\[
\begin{split}
    QG_Y (u) = \exp \Bigg [-\frac{1}{Q_X \left(1 - (1 - u)^{\frac{1}{\theta}}\right)} \int_{0}^{u} &\frac{1}{\theta}q_X \left(1 - (1 - p)^{\frac{1}{\theta}}\right)(1 - p)^{\frac{1}{\theta} - 1}\\& \ln \left(\frac{1}{\theta} q_X \left(1 - (1 - p)^{\frac{1}{\theta}}\right)(1 - p)^{\frac{1}{\theta}}\right)dp\Bigg],
\end{split}
\]
or equivalently
\[
QG_Y (u) = \exp \Bigg[- \frac{1}{Q_X (u)}\int_{0}^{1 - (1 - u)^{\frac{1}{\theta}}}\frac{1}{\theta} q_X (p) (1 - p)^{1 - \theta} \ln \left((1 - p)\frac{1}{\theta}q_X (p)\right)\Bigg] \neq \theta~\mathcal{QG}_X (u).
\]
Also, the quantile version of HFR becomes
\[
QH_Y (u) = \Bigg(\frac{1}{Q_X \left(1 - (1 - u)^{\frac{1}{\theta}}\right)}\int_{0}^{u} \frac{1}{\theta}(1 - p)^{\frac{2}{\theta} -1}q_X \left(1 - (1 - u)^{\frac{1}{\theta}}\right)^2\Bigg)^{-1} \neq \theta ~QH_X (u).
\]
This clearly illustrates that quantile version of the AFR, GFR and HFR for the PHM not satifying the properties those hold in the distribution function approach.
\section{Conclusion}
At the long last, for  readers we reiterate that mixture of $n$ distributions is a special case of formation of $n$ independent component series system having weighted failure rates with the sum of weight functions being unity. However, the latter system having arbitrary weights is also not a generalization of the former. The idea of relating the said concepts deserves some credit because the existing literature on mixture of distributions can be extended to the formation of coherent systems (in particular, series system) so far as non-preservation properties of reliability operations are concerned. One can generate new distributions using  weighted version of arithmetic, geometric and harmonic means of failure rate. Since, the quantile version of means of hazard rate is a special case of weighted means of failure rate, the properties  studied for weighted means is put forth for the prior.
\section*{Acknowledgement}
Subarna Bhattacharjee  would like to thank Odisha State Higher Education Council for providing support to carry out the research project under OURIIP, Odisha, India (Grant No. 22-SF-MT-073).

\end{document}